\documentclass[12pt]{amsart}
\usepackage[a4paper, left=3cm, right=3cm, height=22cm]{geometry}
\usepackage{amsmath,amsfonts,amsthm, amssymb}
\usepackage[arrow,matrix, all,cmtip]{xy}
\usepackage{enumerate}

\usepackage[usenames]{color}
\usepackage{amssymb}

\newtheorem{theorem}{Theorem}[section]
\numberwithin{equation}{section}

\theoremstyle{remark}

\newtheorem{conjecture}[equation]{Conjecture}

\def\Z{{\mathbb{Z}}}

\usepackage[pdfauthor={Jin Cao}, %
				pdftitle={Cdga}%
			dvips,colorlinks=true]{hyperref}
\hypersetup{
	bookmarksnumbered=true,
	linkcolor=black,
}
\newcounter{elno}

\begin{document}
\author{Jin Cao and Wenchuan Hu}
\title{The equivalence of Friedlander-Mazur  and Standard conjectures for threefolds}
\date{\today}

\address{
Jin Cao\\
Yau Mathematical Sciences Center,
Tsinghua University, Beijing, China
}
\email{caojin@mail.tsinghua.edu.cn}

\address{Wenchuan Hu,
School of Mathematics,
Sichuan University,
Chengdu,
China
}
\email{huwenchuan@gmail.com}


\begin{abstract}
 We show that the Friedlander-Mazur conjecture holds for a complex smooth projective variety $X$ of dimension three implies the standard conjectures hold for $X$. This together with a result of Friedlander yields the equivalence of the two conjectures in dimension three.
 From this we provide some new examples whose standard conjectures hold.
\end{abstract}
\maketitle
\pagestyle{myheadings}
 \markright{The equivalence of Friedlander-Mazur  and Standard conjectures for threefolds}

\section{The main result}
Let $X$ be a complex projective variety of dimension $n$. The Lawson homology
$L_pH_k(X)$ of $p$-cycles for a projective variety is defined by
$$L_pH_k(X) := \pi_{k-2p}({\mathcal Z}_p(X)) \quad for\quad k\geq 2p\geq 0,$$
where ${\mathcal Z}_p(X)$ is the space of algebraic $p$-cycles on $X$ provided with a natural topology (see
\cite{Friedlander1}, \cite{Lawson1}).

In \cite{Friedlander-Mazur}, Friedlander and Mazur define a cycle class map
\begin{equation}\label{eq01}
\Phi_{p,k}:L_pH_{k}(X)\rightarrow H_{k}(X,\Z)
\end{equation}
for all $k\geq 2p\geq 0$. Suslin conjectures  that $\Phi_{p,k}:L_pH_{k}(X)\rightarrow H_{k}(X,\Z)$ is an injection for $k=n+p-1$ and an isomorphism for $k\geq n+p$.

From now on, all the (co)homology have the rational coefficients. Let $X$ be a complex projective variety of dimension 3.
In this case, the second author showed in \cite{Hu} that the Friedlander-Mazur conjecture \cite[\S 7]{Friedlander-Mazur} can be reduced to (equivalent to):
\begin{conjecture} \label{FM conj}
\textbf{$\Phi_{1,4}$ is surjective.}
\end{conjecture}

We want to study the relations of this conjecture to  Grothendieck standard conjectures that are fundamental conjectures in the theory of algebraic cycles. We refer to \cite{Kleiman} for the statements of standard conjectures. Recall that the standard conjectures reduce to the standard conjecture of the Lefschetz type.
\

We borrow notations from \cite[Section 2.3]{Beilinson}:
\begin{itemize}
\item
$L(W)$: the Lefschetz type standard conjecture for $W$;
\item
$l(W)$: for every $i > 0$, one can find a  finite correspondence on $W$ that yields an isomorphism $H^{\dim W+ i}(W) \cong H^{\dim W- i}(W)$;
\item
$S(W)$: For any $ j > 0$, one can find a finite correspondence $f_j: W \to Y_j$ with $Y_j$ projective smooth of dimension $j$ such that $f_j^*: H^j(Y_j) \to H^j(W)$ is surjective.
\item
Let $L(n)$ be the assertion that $L(W)$ is true for all $W$ of dimension $\leq n$; same for $l(n), S(n)$.
\end{itemize}
Note that for any smooth projective variety $Y$ with $\dim Y\leq 2$, $L(Y)$ is known. See \cite{Lewis} for example. Therefore $L(2)$ holds.
\

In \cite[Section 2.3]{Beilinson}, Beilinson has proved that $S(n)$ implies $L(W)$ for any smooth projective variety $W$ of dimension $n$ via two steps:
\begin{itemize}
\item
$S(n)$ and $L(n-1)$ implies $l(n)$;
\item
$l(W)$ implies $L(W)$.
\end{itemize}
However the same proof of the first step can be used to prove that $S(X)$, where $X$ is a smooth projective threefold, implies $l(X)$ since $L(2)$ holds. Combining this with his proof of the second step $l(X) \Rightarrow L(X)$, we conclude that $S(X) \Rightarrow L(X)$ for a smooth projective threefold $X$.
\

On the other hand, the Friedlander-Mazur conjecture (\ref{FM conj}) holding for $X$ implies that $S(X)$ via the proof of (i) $\Rightarrow$ (iii) in \cite[Proposition 2.2]{Beilinson}.

In summary, we have the following result.
\begin{theorem}
  The Friedlander-Mazur conjecture for a smooth projective threefold $X$ implies that the Grothendieck standard conjecture holds for $X$.
\end{theorem}

It is already known that Grothendieck standard conjecture holds for a projective threefold  $X$  implies that the Friedlander-Mazur conjecture for $X$ (see \cite[\S 4]{Friedlander2}). Hence  the two conjectures are equivalent in dimension three.

We remark that the Friedlander-Mazur conjecture for all smooth projective varieties is equivalent to that the Grothendieck standard conjecture
of the Lefschetz type  holds.
We emphasis here whether the Friedlander-Mazur conjecture for a given smooth projective variety $X$ is equivalent to the Grothendieck standard conjecture of the Lefschetz type  is still open in dimension four or higher.

\section{Examples}
In most earlier literatures, one verified the Grothendieck standard conjecture first, then use it to  check the Friedlander-Mazur conjecture.
For some  threefolds, it is easier to check the  Friedlander-Mazur conjecture than to check the standard conjectures. This is one of the main purpose in this note.

In \cite{Hu}, it has been proved that the Friedlander-Mazur conjecture holds for a smooth projective threefold $X$ with $h^{2,0}(X) = 0$.
 Hence the standard conjectures also hold for such an $X$. For example,  a smooth projective threefold with representable Chow
 group ${\rm Ch}_0(X)$ of zero cycles satisfying $h^{2,0}(X)=0$. Fano threefolds are examples in this case. A Calabi-Yau threefold $X$ also satisfies $h^{2,0}(X)=0$.
 The last two examples are of threefolds with the Kodaira dimension less than three, the Grothendieck standard conjectures have been proved by
 Tankeev \cite{Tankeev} in a different method.

Surely there exist  lots of threefolds $X$ of general type satisfying $h^{2,0}(X)=0$ but not a complete intersection or a product of a projective curve and surface. Concrete examples of this type can be found in  \cite{Cristina} and references therein. They are new examples of threefolds
for which the Grothendieck standard conjectures of the Lefschetz type hold.
\


\end{document}